\documentclass[]{interact}
\usepackage{latexsym}
\usepackage{anyfontsize} 
\usepackage{tikz} 
\usetikzlibrary{cd} 
\usepackage{color}
\usepackage{hyperref}
\usepackage{url}
\usepackage{breakurl}
\newcommand{\bburl}[1]{\textcolor{blue}{\url{#1}}}

\makeatletter
\newcommand{\monthyear}[1]{%
  \def\@monthyear{\uppercase{#1}}}
\newcommand{\volnumber}[1]{%
  \def\@volnumber{\uppercase{#1}}}
\makeatother

\theoremstyle{plain}
\numberwithin{equation}{section} 
\newtheorem{thm}{Theorem}[section]
\newtheorem{theorem}[thm]{Theorem}

\newtheorem*{theorem*}{Theorem}

\numberwithin{table}{section} 
\numberwithin{figure}{section}

\usepackage{cleveref}

\newcommand{\Id}{\mathrm{Id}}

\newcommand{\av}[1]{\left|#1\right|}
\newcommand{\set}[1]{\left\{#1\right\}}

\definecolor{darkolivegreen}{rgb}{0.33, 0.42, 0.18}

\newcommand\rquot[2]{
  \mathchoice
  {
    \text{\raise0.5ex\hbox{$#1$}\big/\lower0.5ex\hbox{$#2$}}%
  }
  {
    #1\,/\,#2
  }
  {
    #1\,/\,#2
  }
  {
    #1\,/\,#2
  }
}

\newcommand\lquot[2]{
  \mathchoice
  {
    \text{\lower0.5ex\hbox{$#1$}\big\backslash\raise0.5ex\hbox{$#2$}}%
  }
  {
    #1\,\backslash\,#2
  }
  {
    #1\,\backslash\,#2
  }
  {
    #1\,\backslash\,#2
  }
}

\newcommand\lrquot[3]{
  \mathchoice
  {
    \text{\lower0.5ex\hbox{$#1$}\big\backslash\raise0.5ex\hbox{$#2$\!}\big/
      \lower0.5ex\hbox{\!\!$#3$}}%
  }
  {
    #1\,\backslash\,#2\,/\,#3
  }
  {
    #1\,\backslash\,#2\,/\,#3
  }
  {
    #1\,\backslash\,#2\,/\,#3
  }
}

\newcommand{\F}{\mathbb{F}}

\newcommand{\N}{\mathbb{N}}

\newcommand{\Z}{\mathbb{Z}}

\theoremstyle{plain}

\newtheorem*{claim*}{Claim}

\newcommand{\smat}[1]{\left( \begin{smallmatrix} #1 \end{smallmatrix} \right)}
\newcommand{\Fib}{\mathrm{Fib}}
\newcommand{\PhiOrder}{m}
\newcommand{\PeriodOrder}{\frac{p-1}{m}}

\begin{document}


\title{Fibonacci sequences in $\F_p$}

\author{
\name{Menny Aka\thanks{Email: menny.akka@math.ethz.ch}}
\affil{Department of Mathematics, ETH Z\"urich, R\"amistrasse 101, 8092 Z\"urich, Switzerland}
}


\thanks{This research was supported by the SNF Grant 10003145, Equidistribution in Number Theory.}
\maketitle

\begin{abstract}
We prove a conjecture due to Gica \cite{Gica}. Our proof is simple, uses only elementary abstract algebra, and generalizes to other recursive sequences.
\end{abstract}


\section{Introduction}

Let $p$ be a prime number and $\F_p$ (resp.~$\F_{p^2}$) be the finite field of order $p$ (resp.~$p^2$). 
Any choice of $a_1, a_2 \in \F_p$ gives rise to a sequence $F_{a_1,a_2} = (a_k)_{k=1}^\infty$ defined by
\[
A^k \begin{pmatrix} a_1 \\ a_{2} \end{pmatrix} = \begin{pmatrix} a_{k+1} \\ a_{k+2} \end{pmatrix}, \quad \forall k \in \N,
\]
where $A = \smat{0 & 1 \\ 1 & 1}$. Let $\Fib(\F_p)$ denote the set of all such sequences, and $\Fib^*(\F_p)$ the set of sequences $(a_n)$ with $a_n \neq 0$ for all $n \in \N$. As $A$ has finite multiplicative order modulo $p$, any sequence in $\Fib(\F_p)$ is periodic.

Let $\varphi, \varphi'$ denote the eigenvalues of $A$, i.e., the roots of $x^2-x-1=0$, and note that $\varphi'=-\varphi^{-1}$. As this polynomial has discriminant $5$, it follows (e.g.~from quadratic reciprocity) that for $p \neq 2,5$,
\begin{equation}
\varphi, \varphi' \in \F_{p} \quad \Longleftrightarrow \quad \left( \frac{p}{5} \right) = 1 \quad \Longleftrightarrow \quad p \equiv 1,4 \pmod{5},
\label{eq:Recip}
\end{equation}
while $\varphi, \varphi' \in \F_{p^2}\setminus\F_p$ if and only if $\left( \frac{p}{5} \right) = -1$, i.e.~$p \equiv 2,3 \pmod{5}$.


The aim of this note is to resolve a conjecture phrased in \cite{Gica}. As our proof uses only general algebraic arguments, it allows for many generalizations; we propose some extensions in Section~\ref{sec:generalizations}. We denote the order of an element $g$ in group $G$ by $\operatorname{ord}_{G}(g)$ and set $\left( \F_p^\times \right)^{r} = \left\{ a^{r} \;\middle|\; a \in \F_p^\times \right\}$.

\begin{theorem}\label{thm:main}
Let $p\neq 2,5$ be a prime number and $\PhiOrder\in \N$ with $\PhiOrder\mid p-1$. The following are equivalent.
\begin{enumerate}
    \item\label{item:PowerGroup} There exists $(a_n)_{n=0}^\infty \in \Fib^*(\F_p)$ with $\{ a_n \mid n \in \N \} = \left( \F_p^\times \right)^{\PeriodOrder}$.
    
    \item\label{item:Periodm} There exists $(a_n)_{n=0}^\infty \in \Fib^*(\F_p)$ with $\PhiOrder$ as its minimal period length.

    \item \label{item:SplitAndOrder} The prime $p$ satisfies $p \equiv 1,4 \pmod{5}$, and $\operatorname{ord}_{\F_p^\times}(\varphi) = \PhiOrder \lor \operatorname{ord}_{\F_p^\times}(\varphi') = \PhiOrder$.
    
\end{enumerate}

\end{theorem}

\section{Proof of Theorem \ref{thm:main}}
To help the motivated reader without much background in abstract algebra, we give concrete references below to the books \cite{DummitFoote} and \cite{lidl1997finite}.

$(\ref{item:SplitAndOrder}) \implies (\ref{item:PowerGroup})$:
Assume first that $\operatorname{ord}_{\F_p^\times}(\varphi) = \PhiOrder$. Note that 
$$(a_n) := F_{1,\varphi} = ( \varphi^n )_{n=0}^\infty \in \Fib^*(\F_p).
$$ 
Therefore we have $$\{ a_n \mid n \in \N \} = \langle \varphi \rangle \leq \F_p^\times$$ is the subgroup generated by $\varphi$. 
By assumption, it has order $\PhiOrder$. Since $\F_p^\times$ is cyclic of order $p-1$ (see \cite[Thm.~2.8]{lidl1997finite}), $\left( \F_p^\times \right)^{\frac{p-1}{d}} = \left\{ a^{\frac{p-1}{d}} \mid a \in \F_p^\times \right\}$ is the unique subgroup of $\F_p^\times$ of order $d$, it follows that $\langle \varphi \rangle = \left( \F_p^\times \right)^{\PeriodOrder}$, as needed (see \cite[\S 2.3, Theorem 7,(3)]{DummitFoote}. If $\operatorname{ord}_{\F_p^\times}(\varphi') = \PhiOrder$ the same argument applies with $F_{1,\varphi'}$.

$(\ref{item:PowerGroup}) \implies (\ref{item:Periodm})$:
This is now obvious as $\left| \left( \F_p^\times \right)^{\PeriodOrder} \right| = \PhiOrder$ (see \cite[\S 2.3]{DummitFoote}).

$(\ref{item:Periodm}) \implies (\ref{item:SplitAndOrder})$:
Set $l = \operatorname{ord}_{\F_{p^2}^\times}(\varphi)$, $l' = \operatorname{ord}_{\F_{p^2}^\times}(\varphi')$, and let $M_0 = \min(l, l')$ and $M_1 = \max(l, l')$. If $l=l'$ then $M_0=M_1$. If $l\neq l'$, say $l< l'$, we use that $\varphi' = -\varphi^{-1}$, and get 
$$
(\varphi')^{2l} = (-\varphi^{-1})^{2l}= (\varphi)^{-2l}=1
$$
so $l'\mid 2l$. As we assumed that $l<l'$ we get $l'=2l$ so $M_1=2M_0$. The case $l<l'$ is similar, so in conclusion we have  either $M_1 = M_0$ or $M_1 = 2 M_0$.

By assumption, we have  $F_{a_1,a_2} \in \Fib^*(\F_p)$ with a minimal period of length $m$. We first establish the second assertion, namely, that $m =l$ or $m = l'$. To this end, note that since $A^{M_1}$ is similar to $\smat{\varphi^{M_1} & 0 \\ 0 & \varphi'^{M_1}} = \Id$, we have $A^{M_1} = \Id$. Therefore, $M_1$ is also a period for $F_{a_1,a_2}$, so $m \mid M_1$. 
Moreover, since $F_{a_1,a_2}$ has period $m$, 
\[
A^m \begin{pmatrix} a_1 \\ a_2 \end{pmatrix} = \begin{pmatrix} a_1 \\ a_2 \end{pmatrix} \neq \begin{pmatrix} 0 \\ 0 \end{pmatrix},
\]
so $1$ is an eigenvalue of $A^m$. 
As the eigenvalues of $A^m$ are $\varphi^m, \varphi'^{m}$, at least one of them is equal to $1$. So either $l \mid m$ or $l' \mid m$.
Suppose $l \mid m$. If $l = M_1$, then $m = l$, as we saw before that $m\mid M_1$. If $l = M_0< M_1$, then $M_0 \mid m$, $m \mid M_1 = 2 M_0$, so either $m = M_0 = l$ or $m = 2 M_0 = l'$, as we wanted to show. The case $l' \mid m$ is completely analogous.

It remains to show the first assertion $p \equiv 1,4 \pmod{5}$, which  is equivalent to 
$\varphi \in \F_p$ and to $\varphi' \in \F_p$. Assume to the contrary that
\begin{equation}\label{eq:Shlila}
\varphi,\varphi'\in \F_{p^2}\setminus\F_p.    
\end{equation}
 Then they belong to the subgroup
\begin{equation}\label{eq:DefNpm}
\varphi,\varphi' \in N^{\pm} := \left\{ a \in \F_{p^2}^\times \mid \operatorname{Nr}(a) = \pm 1 \right\},    
\end{equation}
 where $\operatorname{Nr}(a) = a \cdot a^p$ is the norm map from $\F_{p^2}$ to $\F_p$ (see e.g.~\cite[Def.~2.7]{lidl1997finite}). As the norm map $\operatorname{Nr}:\F_{p^2}^\times \to \F_{p}^\times$ is surjective (\cite[Thm.~2.8]{lidl1997finite}) the subgroup $N^\pm$ is of size
\[
2\frac{\av{\F_{p^2}^\times}}{\av{\F_{p}^\times}}=2\frac{p^2 - 1}{p - 1} = 2(p + 1).
\]

We consider now three cases: $m=M_0<M_1$, $m=M_0=M_1$, and $M_0<M_1=m$. Starting with  $\PhiOrder=M_0<M_1$, we let $b\in\set{\varphi,\varphi'}$ the element with $\operatorname{ord}_{\F_{p^2}^\times}(b)=2m=M_1$. Since $b\in N^\pm$, $2m\mid2(p+1)$, and together with $m\mid p-1$ we get  $2m\mid 2(p+1)-2(p-1)=4$, so $m=1$ or $m=2$. This is a contradiction to \eqref{eq:Shlila}, as $\pm 1$ are the only elements of order $1$ and $2$ and they belong to $\F_p$.

Now consider the case $M_1=M_2$. We claim that here we must have $m\equiv 0 \pmod 4$. Indeed, if not, either $l=l'=m\equiv 1,3\pmod 4$ or $l=l'=m\equiv 2\pmod 4$. In the former, we have $1=\varphi^l=(-1)^l(\varphi')^{-l'}=-1$, which is absurd. In the latter, we have that $\varphi^{l/2},\varphi'^{l'/2}$ both have order 2, so $$-1=\varphi^{l/2}=(-1)^{l/2}(\varphi')^{-l'/2}=(-1)^2=1,$$ which is also absurd. This shows the claim and therefore $4\mid m\mid p-1$, so $p\equiv 1 \pmod 4$. Similarly to the first case, if we assume that $\varphi,\varphi'\notin \F_p$ we have that $m$ divides both $2(p+1)$ and $p-1$, so it divides $2(p+1)-2(p-1)=4$. If $m=1,2$ we get a contradiction as above. If $m=4$, then $\varphi^4=1$, that is, $\varphi$ is a square root of $-1$. But since $p\equiv 1 \pmod 4$, by quadratic reciprocity $\left( \frac{-1}{p} \right) = 1$, so $\varphi\in \F_p$.

Finally, consider the case $M_0<M_1=m$. Here, $m$ must be even. Assume first that $m\equiv 2\pmod 4$ and suppose that $l=m$ (the case $l'=m$ is analogous). The element $\varphi^2$ has order $\frac{l}{2}=\frac{m}{2}$ and belongs to the group $N^+:=\left\{ a \in \F_{p^2}^\times \mid \operatorname{Nr}(a) = 1 \right\}<N^\pm$, which has size $p+1$. Therefore $\frac{m}{2}\mid p+1$ and since $\frac{m}{2}$ is odd and $p\neq 2$, $m\mid p+1$. It follows that $m\mid (p+1)-(p-1)=2$, which yields a contradiction as above. If $m\equiv 0 \pmod 4$ and say $l=M_0=\frac{m}{2}, l'=m=M_1$, then $1=\varphi^l=(-1)^l(\varphi')^{-l}=(\varphi')^{-\frac{l'}{2}}=-1$, a contradiction.

We have shown that in all possible cases we have $\varphi,\varphi' \in \F_p^\times$, or equivalently that $p \equiv 1,4 \pmod{5}$, which concludes the proof. \qed

\section{Remarks and generalizations} \label{sec:generalizations}

We conclude with several remarks and suggestions for generalizations.  

Theorem \ref{thm:main} is slightly stronger than the conjecture stated in \cite{Gica}, which only asks about the equivalence of items (\ref{item:PowerGroup}) and (\ref{item:SplitAndOrder}). 

The reader may wonder what happens when $p=2$ or $p=5$. Note that for any fixed prime $p$ the set $\Fib^*(\F_p)$ is finite, so one can study it concretely. For $p=2$, it follows readily that $\Fib^*(\F_2)=\emptyset$. For $p=5$ one first checks directly that for $a\in \F_5$ we have $F_{1,a}\in \Fib^*(\F_5)$ if and only if $a=3$. It follows that 
$$
\Fib^*(\F_5)=\set{F_{\alpha,3\alpha}:\alpha\in \F_5^\times}=\set{F_{1,3},F_{3,4},F_{4,2},F_{2,1}},
$$
which are all the shifts of the periodic sequence $\overline{1342}$.

The proof above can be adapted to proving analogous statements for other recursive sequences  in $\F_p$ of the form
\[
a_n = P a_{n-1} - Q a_{n-2},\quad a_1,a_2\in \F_p
\]
with $P^2 - 4Q > 0$ and $Q = \pm 1$ (or more generally, $(P,Q)=1$), that is to Lucas sequences over $\F_p$. For this one studies the order of 
\[
B_{P,Q} = \begin{pmatrix} 0 & 1 \\ -Q & P \end{pmatrix}
\]
 which plays the same role as $A= \smat{0 & 1 \\ 1 & 1}$ for these generalized second-order recursions. In fact, the proof simplifies when $Q=1$, as in this case the eigenvalues $\lambda, \lambda'$ of $B_{P,Q}$ are inverses of each other (so $\operatorname{ord}_{\F_{p^2}^\times}(\lambda) = \operatorname{ord}_{\F_{p^2}^\times}(\lambda')$ for every $p$).

Moreover, the same methods can prove a complementary result to Theorem \ref{thm:main} (or to an analog of it for recursion sequences as above). We explicate it here for the Fibonacci sequence. Recall the definition of $N^{\pm}$ in \eqref{eq:DefNpm}.
\begin{theorem*}
For a prime $p \neq 2,5$ and $m \mid 2(p+1)$, the following are equivalent.
\begin{enumerate}
    \item There exists $(a_n)_{n=1}^\infty \in \Fib^*(\F_p)$ with $\{ a_n \mid n \in \N \} = \left( N^{\pm} \right)^{\frac{2(p+1)}{m}}$.
    
    \item There exists $(a_n)_{n=1}^\infty \in \Fib^*(\F_p)$ with minimal period length $m$.
    
    \item The prime $p$ satisfies $p \equiv 2,3 \pmod{5}$ and  $m = \operatorname{ord}_{\F_{p^2}^\times}(\varphi)\lor m = \operatorname{ord}_{\F_{p^2}^\times}(\varphi')$.
\end{enumerate}
\end{theorem*}

Finally, we remark that this proof also allows considering general moduli $N$ instead of $p$, that is, one can formulate and prove similar results working with $(\Z/N\Z)^\times$ instead of $\F_p^\times=(\Z/p\Z)^\times$.

While these generalizations can be proven with the same ideas, they require significantly more technical detail. We leave them as challenges to the motivated reader.

\bibliographystyle{amsalpha}
\bibliography{Bibliography}

\providecommand{\bysame}{\leavevmode\hbox to3em{\hrulefill}\thinspace}
\providecommand{\MR}{\relax\ifhmode\unskip\space\fi MR }
\providecommand{\MRhref}[2]{%
  \href{http://www.ams.org/mathscinet-getitem?mr=#1}{#2}
}
\providecommand{\href}[2]{#2}
\begin{thebibliography}{Gic08}

\bibitem[DF04]{DummitFoote}
David~Steven Dummit and Richard~M Foote, \emph{Abstract algebra}, vol.~3, Wiley Hoboken, 2004.

\bibitem[Gic08]{Gica}
Alexandru Gica, \emph{Quadratic residues in fibonacci sequences}, The Fibonacci Quarterly \textbf{46} (2008), no.~1, 68--72.

\bibitem[LN97]{lidl1997finite}
Rudolf Lidl and Harald Niederreiter, \emph{Finite fields}, no.~20, Cambridge university press, 1997.

\end{thebibliography}
\end{document}